\documentclass[11pt]{amsart}
\usepackage{a4, amsmath,  amsthm,  amsfonts, mathrsfs, latexsym, amssymb, 
pstricks, pst-grad, graphicx, tabularx, setspace, longtable}
\usepackage{enumerate}
\usepackage{hyperref}
\usepackage{color}
\usepackage{psfrag}
\usepackage{amssymb}
\scshape
\date{}
\usepackage{tikz}
\usetikzlibrary{arrows}
\usepackage{subcaption}
\usetikzlibrary{automata}
\newtheorem{thm}{\bf Theorem}[section]
\newtheorem{cor}[thm]{\bf Corollary}

\newtheorem{defn}[thm]{\bf Definition}

\newtheorem{rem}[thm]{\bf Remark}

\begin{document}
\title[Laplacian Eigen values of  character degree graphs of solvable groups]{LAPLACIAN EIGEN VALUES OF  CHARACTER DEGREE GRAPHS OF SOLVABLE GROUPS}
 
\author [G. Sivanesan]{G. Sivanesan}
\address{Department of Mathematics, Government College of Engineering, Salem 636011, Tamil Nadu, India, ORCID: 0000-0001-7153-960X.}
\email{sivanesan@gcesalem.edu.in}

\author [C. Selvaraj]{C. Selvaraj}
\address{Department of Mathematics, Periyar University, Salem 636011, Tamil Nadu, India,  ORCID:  0000-0002-4050-3177.}
\email{selvavlr@yahoo.com}

 
\begin{abstract} Let $G$  be a finite solvable group,  let $Irr(G)$ be the set of all complex irreducible characters of $G$ and let $cd(G)$ be the set of all degrees of characters in $Irr(G).$  Let $\rho(G)$ be the set of primes that divide degrees in $cd(G).$ The character degree graph $\Delta(G)$ of $G$ is the simple undirected graph with vertex set $\rho(G)$ and in which two distinct vertices $p$ and $q$ are adjacent  if there exists a  character degree $r \in cd(G)$ such that $r$ is divisible by the product $pq.$ In this paper, we obtain Laplacian eigen values and distance Laplacian eigen values of regular character degree graph, super graphs of regular character degree graph and  character degree graph with diameter $2$ has two blocks.
\end{abstract}
\date{}
\maketitle
 
\section{Introduction} In this paper, $G$ is a finite solvable group with identity $1$. The set of complex irreducible characters of $G$ is denoted by $Irr(G)$. $cd(G) =\{\chi(1) \mid \chi\in Irr(G)\}$ is the set of all distinct degrees of irreducible characters in $Irr(G)$. Let $\rho(G)$  be the set of all primes that divide degrees in $cd(G)$. An extensive literature has been devoted to studying how graphs can be associated with groups, and that literature can be used to investigate the algebraic structure of groups using graph theoretical properties. It is one of these graphs that is called the character degree graph $\Delta(G)$ of $G$. In fact,  $\Delta(G)$ is an undirected simple graph with vertex set $\rho(G)$ in which $p,q\in\rho(G)$ are joined by an edge if there exists a character degree $\chi(1)\in cd(G)$ which is divisible by $pq.$ It was first defined in  \cite{Manz2} and studied by many authors  (see \cite{PP}, \cite{Zha}, \cite{Lewis3}). There have been some interesting results obtained on the character graph of $G$ when $G$ is a solvable group. Actually, Manz \cite{Manz1}  proved that $\Delta(G)$  has at most two connected components. Also, Manz et al. \cite{Manz3} have proved that diameter of $\Delta(G)$ is at most $3$.  We obtained a necessary condition for the character degree graph $\Delta(G)$  of a finite solvable group $G$ to be Eulerian  \cite{SST}.  Mahdi Ebrahimi et al. \cite{EI} studied all regular character-graphs whose eigenvalues are in the interval $ [-2,\infty)$. Furthermore, Hafezieh et al. studied the effects of cut vertices and eigenvalues of $\Delta(G)$ on the group structure of G \cite{HH}. This study on character degree graphs of solvable groups motivated us to study Laplacian eigenvalues and distance Laplacian eigenvalues of regular character degree graphs, super graphs of regular character degree graphs, and character degree graphs with diameter $2$ has two blocks.

\section{Preliminaries}
This section presents some preliminary results that are used in the paper. The graphs are assumed to be simple, undirected, and finite. Let $\Gamma$ be a graph with vertex set  $V(\Gamma)$ and edge set $E(\Gamma)$. If $\Gamma$ is connected, then the distance $d(u, v)$ between two distinct vertices $u,v\in V(\Gamma)$ is the length of the shortest path between them. An graph's diameter is the maximum distance between all pairs of distinct vertices. A graph with $n$ vertices in which any two distinct vertices are adjacent is denoted by $ K_n$.  $C_n$denotes a cycle with n vertices. The vertex connectivity $k(\Gamma)$ of $\Gamma$ is defined to be the minimum number of vertices whose removal from $\Gamma$ results in a disconnected subgraph of $\Gamma.$ A cut vertex v of a graph $\Gamma$ is a vertex such that the number of connected components of $\Gamma -v$  is more than the number of connected components of $\Gamma$.  A block is a maximally connected subgraph that does not have a cut vertex. By their maximality, different blocks of $\Gamma$ overlap in at most one vertex, which is then a cut vertex.  A vertex's degree is defined as the number of edges incident with it, and it is denoted by $ d(v)$  or  $deg v$. A graph $\Gamma$ is called $k$-regular, if the degree of each vertex is $k$. Cliques of graphs are complete subgraphs and the number of vertices on a clique of maximum size is called its clique number and it is denoted by $\omega(\Gamma)$.  The spanning tree T of an undirected graph G is a subgraph that includes all of its vertices.
\newline For a finite simple undirected graph $\Gamma$ with vertex set $V(\Gamma) =  \{v_1, v_2, . . . , v_n\}$, the adjacency matrix  $A(\Gamma )$ is the $n \times n$ matrix with $(i , j )$th entry is 1 if  $v_i$ and $ v_ j$ are adjacent and $0$ otherwise. We denote the diagonal matrix $D(\Gamma ) = diag(d_1, d_2, . . . , d_n)$  where $d_i$ is the degree of the vertex $v_i$ of $\Gamma$ , $i = 1, 2, . . . , n$. The Laplacian matrix $L(\Gamma)$ of  $\Gamma$ is the matrix $D(\Gamma ) - A(\Gamma).$  The matrix $L(\Gamma)$ is symmetric and positive semidefinite, so that its eigenvalues are real and nonnegative. Furthermore, the sum of each row (column) of  $L(\Gamma)$ is zero. The eigenvalues of $L(\Gamma )$ are called the Laplacian eigenvalues of  $\Gamma$  and it is denoted by  $\lambda_1(\Gamma ) \geq \lambda_2(\Gamma ) \geq · · ·\geq \lambda_n(\Gamma) = 0$. Now let $\lambda_{n_1}(\Gamma)  \geq \lambda_{n_2}(\Gamma)  \geq · · ·  \geq \lambda_{n_r} (\Gamma ) = 0$ be the distinct eigenvalues of $\Gamma$ with multiplicities
$m_1,m_2, . . . ,m_r$, respectively. The Laplacian spectrum of $\Gamma$, that is, the spectrum of $L(\Gamma )$  is represented as
$ \begin{pmatrix}
\lambda_{n_1}(\Gamma) &\lambda_{n_2}(\Gamma) & \cdots & \lambda_{n_r} (\Gamma ) \\
m_1 & m_2 & \cdots & m_r
\end{pmatrix} $.

The distance between two vertices $u, v \in V $ is denoted by  $d_{uv}$ and is defined as the length of a shortest path between $u$ and $v$ in $G$. The distance matrix of $G$ is denoted by $D(G)$ and defined by $D(G) = (d_{uv})$   $u, v \in V $.  The transmission $Tr(v)$ of a vertex $v$ is defined to be the sum of the distances from $v$ to all other vertices in $G$, i.e.

\hspace {5cm} $Tr(v) = \sum_{u \in V} d_{uv} $.
\newline The distance Laplacian of a connected graph $G$ is  $D^{L} (G) =Tr (G) - D (G),$  where $Tr(G)$ denotes the diagonal matrix with  $Tr(v_i)$ as the $i$-th diagonal entry.

We will take into account the following well known facts concerning character degree graphs and they are needed in the next section.

\begin{rem}\normalfont\label{rem2.1} For results regarding $\Delta(G),$ we start with P\'{a}lfy's three prime theorem on the character degree graph of solvable groups. P\'{a}lfy theorem~\cite[Theorem, p. 62]{PP} states that given a solvable group $G$ and any three distinct vertices of $\Delta(G)$, there exists an edge incident with the other two vertices. On applying P\'{a}lfy's theorem, $\Delta(G)$  has at most two connected components. 
\end{rem}

\begin{rem}\normalfont\label{rem2.3} 
Huppert  \cite[p. 25]{HU} listed all possible graphs $\Delta(G)$ for solvable groups $G$ with at most $4$ vertices.   In fact, every graph with $3$ or few vertices that satisfies P\'{a}lfy's condition occurs as $\Delta(G)$ for some solvable group $G$  \cite[p. 184]{Lewis3}. 
\end{rem}
\begin{thm}\cite[Lemma 2.7]{EI}\label{2.2} Let $G$ be a group with $\mid\rho(G)\mid\geq 3.$  If $\Delta(G)$ is not a block and the diameter of $\Delta(G)$ is at most $2$, then each block of $\Delta(G)$ is a complete graph.
\end{thm}

\begin{thm} \cite[Theorem 5]{Zha}\label{2.3}  The graph with four vertices in Figure 1 is not the character degree graph of a solvable group.
	
\begin{figure}[ht]
\centering
\begin{tikzpicture}
\draw[fill=black] (0,0) circle (3pt);
 \draw[fill=black] (1,0) circle (3pt);
 \draw[fill=black] (2,0) circle (3pt);
 \draw[fill=black] (3,0) circle (3pt);
\draw[thick] (0,0) --(1,0) --(2,0) --(3,0);
\end{tikzpicture}
 \caption{Graph with four vertices}
\end{figure}
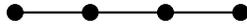
\end{thm}

\begin{thm}\cite[Theorem 1.1]{Lewis4}\label{2.4}  Let $G$ be a solvable group. Then $\Delta(G)$ has at most one cut vertex.
\end{thm}	

\begin{thm}\cite[Theorem A]{MZ}\label{2.5} If $\Delta(G)$ is a non-complete and regular character degree graph of a finite solvable group $G$ with $n$ vertices, then $\Gamma(G)$ is $n-2$ is a regular graph.
\end{thm}

\begin{thm} \cite [Lemma 2.1]{Lewis2}\label{2.6}  Let $p, q_1, . . . , q_n$ be distinct primes so that $p$ is odd. Then there is a
solvable group $G$ of fitting height $2$ so that $\Delta(G)$ has two connected components: $\{p\}$ and $\{q_1, . . . , q_n\}$.
\end{thm}

\begin{thm} \cite [ Theorem A] {Lewis2}\label{2.7} Let  $\Gamma$ be a graph with $n$ vertices. There exists a solvable group $G$ of fitting
height  $2$  with  $\Delta(G) = \Gamma $ if and only if the vertices of degree less than $n - 1$ can be partitioned into two subsets, each of which induces a complete subgraph of     $\Gamma$  and one of which contains only vertices of degree $n - 2$.

\end{thm}

\begin{thm} \cite [ Corollary 4.2]{Bra}\label{2.8} The number  of spanning trees of the graph G of order $n$  $ = \frac{1}{n}   \lambda_2( G )  \lambda_3( G ) \ldots \lambda_n( G ).$  Here $ \lambda_2( G ),  \lambda_3( G ), \ldots, \lambda_n( G ) $ are Laplacian eigen values of G.

\end{thm} 

\begin{thm}   \cite [ Theorem 3.1] {AH}\label{2.9} Let $G$ be a connected graph on n vertices with diameter $D \leq 2.$ Let  $ \lambda_1 ^{L} \geq  \lambda_2 ^{L} \geq  ...  \geq  \lambda_{n-1} ^{L} >  \lambda_{n} ^{L} = 0 $  be the Laplacian spectrum of $G.$ Then the distance Laplacian spectrum of $G$  is $ 2n -\lambda_{n-1} ^{L} \geq  2n -\lambda_{n-2} ^{L} \geq ... \geq  2n -\lambda_{1} ^{L} > \delta_n ^{L} = 0. $ Moreover, for every  $ i  \in \{ 1, 2, ... , n-1 \}$ the eigenspaces corresponding to $  \lambda_{i} ^{L} $ and to $ 2n -\lambda_{n-i} ^{L} $ are the same. 
\end{thm} 

\begin{rem}{\bf Operation D:}\label{rem2.6} \normalfont Consider a character degree graph $\Delta(G).$ Select two distinct primes $p$ and $q$ such that $p$ is odd prime and $p,q\not\in \rho(G).$ By Theorem~\ref{2.6}, there exists a solvable group $H$ of fitting height $2$ and $\Delta(H)$ is the graph containing two isolated vertices. Now one can make a obtain the direct product $\Delta(G\times H).$  Given $\Delta(G),$ the construction of  such a direct product $\Delta(G\times H)$ is used in the proof of Theorem~\ref{3.1} repeatedly. We denote this process as Operation D.
\end{rem}

\begin{defn} \cite[p. 502]{BL}\label{def2.8} Using direct product, one can find bigger groups from smaller groups. The same may be used to construct higher order character degree graphs. For two groups $A$ and $B$ where $\rho(A)$ and  $\rho(B)$ are disjoint,  we have that $\rho(A\times B) =  \rho(A) \cup \rho(B)$.  Define an edge between vertices $p$ and $q$ in $\rho(A\times B)$  if any of the following is satisfied:
\begin{itemize}
\item [\rm (i)] $p , q\in \rho (A)$ and there is an edge between $p$ and $q$ in $\Delta(A);$
\item [\rm (ii)] $p, q \in \rho (B)$ and there is an edge between $p$ and $q$ in $\Delta(B);$
\item [\rm (iii)] $p\in\rho (A)$ and $q\in\rho (B);$
\item [\rm (iv)]  $p \in\rho (B) $ and $q\in\rho (A).$ 
\end{itemize}
Now we get a higher order character degree graph and it is called direct product.
\end{defn}

\section{Laplacian Eigen values of Character Degree Graphs}
If $\Delta(G)$ is a character degree graph and  $u , v$  are two non adjacent vertices  of $\Delta(G)$,  then $\Delta(G)+uv$ may or may not be character degree graph.
\newline Consider the following example:
\begin{figure}[ht]
    \centering 
     \begin{subfigure}{.5\textwidth} 
         \centering 
         \begin{tikzpicture} 
           \draw[fill=black] (1,0) circle (2pt); 
\draw[fill=black] (2,0) circle (2pt); 
\draw[fill=black] (-1,0) circle (2pt);
 \draw[fill=black] (-2,0) circle (2pt);
\draw[fill=black] (0,1.5) circle (2pt); 
\draw[fill=black] (0,-1.5) circle (2pt); 
\draw[thick] (0,1.5) --(1,0)--(0,-1.5)--(-1,0)--(0,1.5);  
\draw[thick] (0,1.5) --(2,0)--(0,-1.5)--(-2,0)--(0,1.5);  
\draw[thick] (-2,0) --(-1,0);  
\draw[thick] (1,0) --(2,0);  
\draw[thick] (0,1.5) --(0,-1.5);
         \end{tikzpicture} 
         \caption{Character degree graph with 6 vertices} 
     \end{subfigure}
     \begin{subfigure}{.5\textwidth} 
         \centering 
         \begin{tikzpicture} 
             \draw[fill=black] (1,0) circle (2pt); 
\draw[fill=black] (2,0) circle (2pt); 
\draw[fill=black] (-1,0) circle (2pt);
 \draw[fill=black] (-2,0) circle (2pt);
\draw[fill=black] (0,1.5) circle (2pt); 
\draw[fill=black] (0,-1.5) circle (2pt); 
\draw[thick] (0,1.5) --(1,0)--(0,-1.5)--(-1,0)--(0,1.5);  
\draw[thick] (0,1.5) --(2,0)--(0,-1.5)--(-2,0)--(0,1.5);  
\draw[thick] (-2,0) --(-1,0);  
\draw[thick] (1,0) --(2,0);  
\draw[thick] (0,1.5) --(0,-1.5);
\draw[thick] (-1,0) --(1,0);  
         \end{tikzpicture} 
         \caption{Unknown six-vertex graph} 
\end{subfigure} 
\caption{Character degree graph with addition of line}
\end{figure}
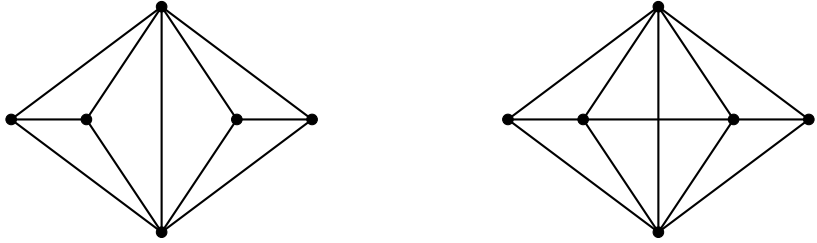 
\newline  Figure  2(A) shows a character degree graph as shown in \cite[p. 502-503]{BL}.  Figure 2(B) cannot be concluded as a character degree graph at this time, as shown in \cite[p. 509]{BL}.
\newline Interestingly, all  graphs which are obtained by adding two non-adjacent vertices in a $n-2$ regular character degree graph  are character degree graphs.

\begin{thm}\label{3.1}
 Let $n\geq 4$ be an even integer. Let $\Delta(G)$ be the  $n-2$ regular character degree graph, then all the super graphs of $\Delta(G)$ by addition of two non adjacent vertices are  character degree graphs  for some solvable group $G.$

\end{thm}
\begin{proof}
Using the induction method, let us prove the result based on the number of vertices in $\Delta(G).$

\textbf{Case 1.} Let $n=4.$ 
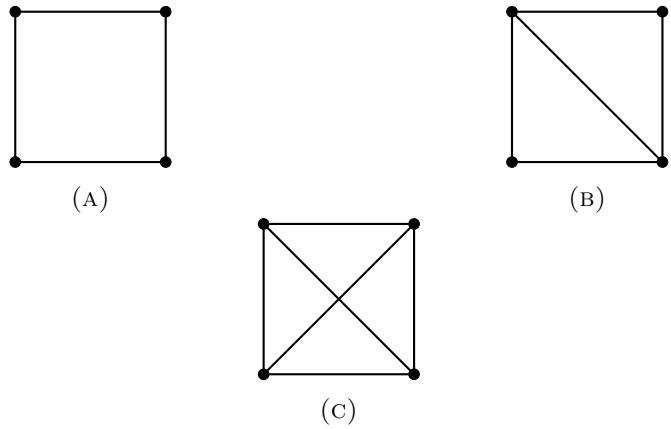
\begin{figure}[ht]
    \centering 
     \begin{subfigure}{.5\textwidth} 
         \centering 
         \begin{tikzpicture} 
           \draw[fill=black] (0,0) circle (2pt); 
 \draw[fill=black] (2,0) circle (2pt); 
 \draw[fill=black] (0,2) circle (2pt); 
 \draw[fill=black] (2,2) circle (2pt); 
\draw[thick] (0,0) --(2,0)--(2,2)--(0,2)--(0,0);
   
         \end{tikzpicture} 
         \caption{} 
     \end{subfigure}
     \begin{subfigure}{.5\textwidth} 
         \centering 
         \begin{tikzpicture} 
            \draw[fill=black] (0,0) circle (2pt); 
 \draw[fill=black] (2,0) circle (2pt); 
 \draw[fill=black] (0,2) circle (2pt); 
 \draw[fill=black] (2,2) circle (2pt); 
\draw[thick] (0,0) --(2,0)--(2,2)--(0,2)--(0,0);
\draw[thick] (0,2) --(2,0);
         \end{tikzpicture} 
         \caption{} 
     \end{subfigure} 
\begin{subfigure}{.5\textwidth} 
         \centering 
         \begin{tikzpicture} 
 \draw[fill=black] (0,0) circle (2pt); 
 \draw[fill=black] (2,0) circle (2pt); 
 \draw[fill=black] (0,2) circle (2pt); 
 \draw[fill=black] (2,2) circle (2pt); 
\draw[thick] (0,0) --(2,0)--(2,2)--(0,2)--(0,0);
\draw[thick] (0,2) --(2,0);
\draw[thick] (0,0) --(2,2);
         \end{tikzpicture} 
         \caption{  } 
     \end{subfigure}
\caption{ Character degree graphs  with four vertices} 
\end{figure} 

Figure $3$ shows character degree graphs as shown in \cite[p. 184]{Lewis3}.  Figures 3(B) and 3(C) are supergraphs of Figure 3(A). By adding two non-adjacent vertices to a 2-regular graph with four vertices, super graphs are formed.

\textbf{Case 2.} Let $n=6.$ 

\begin{figure*}[ht]
    \centering 
     \begin{subfigure}[t]{.4\textwidth} 
         \centering 
         \begin{tikzpicture} 
           \draw[fill=black] (2,3) circle (2pt); 
\draw[fill=black] (0.5,2) circle (2pt);
 \draw[fill=black] (0.5,0) circle (2pt);
 \draw[fill=black] (2,-1) circle (2pt);
 \draw[fill=black] (3.5,0) circle (2pt);
 \draw[fill=black] (3.5,2) circle (2pt);
\draw[thick] (2,3) --(0.5,2)--(0.5,0)--(2,-1)--(3.5,0)--(3.5,2)--(2,3);
\draw[thick] (2,3) --(0.5,0)--(3.5,0)--(2,3);
\draw[thick] (0.5,2) --(2,-1)--(3.5,2)--(0.5,2);   
         \end{tikzpicture} 
         \caption{} 
     \end{subfigure}
     \begin{subfigure}[t]{.4\textwidth} 
         \centering 
         \begin{tikzpicture} 
            \draw[fill=black] (2,3) circle (2pt); 
\draw[fill=black] (0.5,2) circle (2pt);
 \draw[fill=black] (0.5,0) circle (2pt);
 \draw[fill=black] (2,-1) circle (2pt);
 \draw[fill=black] (3.5,0) circle (2pt);
 \draw[fill=black] (3.5,2) circle (2pt);
\draw[thick] (2,3) --(0.5,2)--(0.5,0)--(2,-1)--(3.5,0)--(3.5,2)--(2,3);
\draw[thick] (2,3) --(0.5,0)--(3.5,0)--(2,3);
\draw[thick] (0.5,2) --(2,-1)--(3.5,2)--(0.5,2);
\draw[thick] (2,3) --(2,-1);
         \end{tikzpicture} 
         \caption{} 
     \end{subfigure} 
\begin{subfigure}[t]{.4\textwidth} 
         \centering 
         \begin{tikzpicture} 
 \draw[fill=black] (2,3) circle (2pt); 
\draw[fill=black] (0.5,2) circle (2pt);
 \draw[fill=black] (0.5,0) circle (2pt);
 \draw[fill=black] (2,-1) circle (2pt);
 \draw[fill=black] (3.5,0) circle (2pt);
 \draw[fill=black] (3.5,2) circle (2pt);
\draw[thick] (2,3) --(0.5,2)--(0.5,0)--(2,-1)--(3.5,0)--(3.5,2)--(2,3);
\draw[thick] (2,3) --(0.5,0)--(3.5,0)--(2,3);
\draw[thick] (0.5,2) --(2,-1)--(3.5,2)--(0.5,2);
\draw[thick] (2,3) --(2,-1);
\draw[thick] (0.5,2) --(3.5,0);
         \end{tikzpicture} 
         \caption{  } 
     \end{subfigure}
\begin{subfigure}[t]{.4\textwidth}
         \centering 
         \begin{tikzpicture} 
 \draw[fill=black] (2,3) circle (2pt); 
\draw[fill=black] (0.5,2) circle (2pt);
 \draw[fill=black] (0.5,0) circle (2pt);
 \draw[fill=black] (2,-1) circle (2pt);
 \draw[fill=black] (3.5,0) circle (2pt);
 \draw[fill=black] (3.5,2) circle (2pt);
\draw[thick] (2,3) --(0.5,2)--(0.5,0)--(2,-1)--(3.5,0)--(3.5,2)--(2,3);
\draw[thick] (2,3) --(0.5,0)--(3.5,0)--(2,3);
\draw[thick] (0.5,2) --(2,-1)--(3.5,2)--(0.5,2);
\draw[thick] (2,3) --(2,-1);
\draw[thick] (0.5,2) --(3.5,0);
\draw[thick] (0.5,0) --(3.5,2);
         \end{tikzpicture} 
         \caption{  } 
     \end{subfigure}
\caption{ Character degree graphs with six vertices } 
\end{figure*}
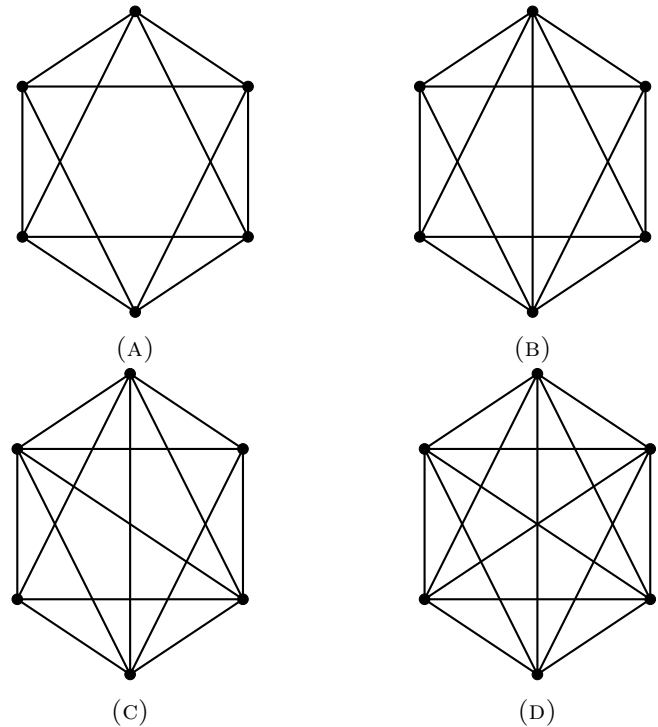 

The graphs in Figure $4$  are character degree graphs as shown in  \cite[p. 502-503]{BL}.  It is important to note that Figures $4(B)$  through Figure $4(D)$  are super graphs of Figure $4(A)$. Supergraphs are created by adding two non-adjacent vertices to four regular graphs with six vertices.

\textbf{Case 3.} Let $n=8.$  Through direct product construction, we can construct a six-regular graph with eight vertices and all of its supergraphs by adding two adjacent vertices.

As in Case $3$, there are four character degree graphs with six vertices.  Let us consider each of these graphs as a $\Delta(G)$.  After applying Operation D~\ref{rem2.6},  we get four direct products $\Delta(G\times H)$  with $8$ vertices.  First, we have a $6$-regular graph with $8$ vertices. In the second graph, two vertices have a degree of $7$ and six vertices have a degree of $6$. On the third graph, four of the vertices are of degree $7$  while the remaining four are of degree $6$. The fourth graph contains six vertices of degree $7$ and two vertices of degree $6$. The complete graph occurs as $\Delta(G)$ for some solvable group G \cite[p. 39]{BL1}.  As a result, a complete graph with $8$ vertices is also a character degree graph. Thus, all the supergraphs of the $6$-regular graph with $8$ vertices are character degree graphs.

\textbf{Case 4.} Let $n=10.$  Four supergraphs of $6$-regular graphs with $8$ vertices and one $6$-regular graph are constructed in case $3$. We can take each of these graphs as a $\Delta(G)$. After applying Operation  D~\ref{rem2.6},  we obtain five direct products  $\Delta(G\times H)$  with $10$ vertices. An 8-regular graph with 10 vertices is the first graph.  There are two vertices of degree $9$ and eight vertices of degree $8$ in the second graph.  The third graph has four vertices of degree $9$ and six vertices of degree $8$. In the fourth graph, there are six vertices of degree $9$ and four of degree $8$. The fifth graph has eight vertices of degree $9$ and two vertices of degree $8$. The complete graph with $10$ vertices is also known as the character degree graph. All supergraphs of the $8$-regular graph with $10$ vertices are character degree graphs.
\newline  Assume that as induction hypothesis, there exists  $\frac{n-2}{2}$ super graphs  of $n-4$ regular graph with $n-2$ vertices and one $n-4$ regular graph with $n-2$ vertices. Let us take each of these graphs as  $\Delta(G)$.  By applying Operation D~\ref{rem2.6}, we obtain $\frac{n}{2}$  direct products  $\Delta(G\times H)$  with $n$ vertices. These $\frac{n}{2}$ direct product graphs are super graphs of $n - 2 $ regular graph with $n$ vertices.  As stated in  \cite[p. 39]{BL1},  the complete graph with $n$ vertices is also a character degree graph. 
\newline Therefore, all super graphs of $n-2$ regular graph with $n$ vertices are character degree graphs.\hfill $\square$

\end{proof}

\begin{rem}\label{rem2.7}  Each vertex in an $n - 2$ regular graph is not adjacent to another vertex. Let us say that the vertex $v_i$ $(1\leq i \leq {\frac{n}{2}} )$ is not adjacent to the vertex $v_{i+{\frac{n}{2}}}$. Let $\Delta(G)$ be the supergraph of a $n-2$ regular graph with $n$ vertices. Initially, we label the vertices of $\Delta(G)$  as $\{ v_1,  v_2,  v_3, ... ,  v_{n-1}, v_n \} $. Assume that $m$ is the number of vertices having a degree of $n - 1$.  Remaining $n - m$ vertices have degree $n-2$  can be partitioned into two subsets. Let  r $ = \frac{n-m}{2}$ . There is no difficulty in finding $r$  vertices that are adjacent to each other. As a result, it forms a complete subgraph of $\Delta(G)$. Let us refer to this vertex set as U $ = \{u_1, . . . , u_r \} $.  Each vertex in $U$ has a degree $n - 2.$  Thus, for each vertex in $U$, there is correspondingly one vertex which is not adjacent to that vertex. We shall label this set of vertex values as $W  = \{w_1, . . . , w_r \} $.  In W, each vertex has a degree of $n-2.$  Therefore, all the vertices in $W$ form a complete subgraph of  $\Delta(G)$. Based on Theorem~\ref{2.7}, G is a solvable group of fitting height $2$  for supergraphs of $n - 2$ regular character degree graph. \hfill $\square$

\end{rem}

The next step is to find the Laplacian eigenvalues and multiplicities of an $n-2$ regular graph.

\begin{thm}\label{3.2}  Let $n\geq 4$ be an even integer and $G$ be a solvable group such that $\Delta(G)$  is  $n-2$ regular graph. The characteristic polynomial of the  Laplacian matrix of $\Delta(G)$  is given by $  x     (x-(n-2))^{\frac{n}{2}}   (x-n)^{{\frac{n}{2}} -1}$ .
\end{thm}
\begin{proof}  The Laplacian matrix of $\Delta(G)$ is  $n\times n$ matrix. Each row and column is indexed based on its vertices  $ v_1,  v_2,  v_3, ... , v_{\frac{n}{2}},  v_{{\frac{n}{2}}+1},  . . . , v_n$. In $n-2$ regular graph,  each vertex is not adjacent with one vertex.  Let's say the vertex  $v_i$ $(1\leq i \leq {\frac{n}{2}} )$ is not adjacent with the vertex  $v_{i+{\frac{n}{2}}}$. Then
    
$L(\Delta(G))= $
  $
	\begin{bmatrix} 
	n-2 & -1 & -1 \cdots -1& 0 &-1&-1  \cdots  & -1  \\
	-1 & n-2 & -1   \cdots  -1 & -1 & 0 &-1 \cdots  & -1 \\
          -1 & -1 &n -2   \cdots  -1 & -1 & -1 &0  \cdots  & -1 \\
           \vdots  & \cdots   & \cdots         & \cdots           & \cdots & \cdots    &   \vdots \\
            0 & -1 & -1 \cdots -1& n-2&-1&-1  \cdots  & -1  \\
            -1& 0 & -1 \cdots -1& -1&n-2&-1  \cdots  & -1 \\
             \vdots  & \cdots       & \cdots           & \cdots & \cdots & \cdots   &   \vdots \\
-1 & -1 & -1 \cdots 0& -1 &-1&-1  \cdots  & n-2 
	\end{bmatrix}
	$
\newline The characteristic polynomial of $L(\Delta(G))$ is  $\mid L(\Delta(G)-xI) \mid$

  $=$  $
	\begin{vmatrix} 
	(n-2)-x & -1 & -1 \cdots -1& 0 &-1&-1  \cdots  & -1  \\
	-1 & (n-2)-x & -1   \cdots  -1 & -1 & 0 &-1 \cdots  & -1 \\
          -1 & -1 & (n-2)-x   \cdots  -1 & -1 & -1 &0  \cdots  & -1 \\
           \vdots  & \cdots   & \cdots         & \cdots           & \cdots & \cdots    &   \vdots \\
            0 & -1 & -1 \cdots -1& (n-2)-x &-1&-1 \cdots  & -1  \\
            -1& 0 & -1 \cdots -1& -1& (n-2)-x &-1 \cdots  & -1 \\
             \vdots  & \cdots   & \cdots              & \cdots & \cdots & \cdots   &   \vdots \\
-1 & -1 & -1 \cdots 0& -1 &-1&-1 \cdots  & (n-2)-x
	\end{vmatrix}
	$

\vspace{0.5cm}

 $=$  $
	\begin{vmatrix} 
	-x & -x &  -x  \cdots  -x  &  -x  & -x  & -x   \cdots  &  -x   \\
	-1 & (n-2)-x & -1   \cdots  -1 & -1 & 0 &-1 \cdots  & -1 \\
          -1 & -1 & (n-2)-x   \cdots  -1 & -1 & -1 &0  \cdots  & -1 \\
           \vdots  & \cdots   & \cdots         & \cdots           & \cdots & \cdots    &   \vdots \\
            0 & -1 & -1 \cdots -1& (n-2)-x &-1&-1 \cdots  & -1  \\
            -1& 0 & -1 \cdots -1& -1& (n-2)-x &-1 \cdots  & -1 \\
             \vdots  & \cdots   & \cdots              & \cdots & \cdots & \cdots   &   \vdots \\
-1 & -1 & -1 \cdots 0& -1 &-1&-1 \cdots  & (n-2)-x
	\end{vmatrix}
	$
 (Apply the row operation  $R_1 \rightarrow R_1 + R_2 + R_3 + ... +R_n$ )
\vspace{0.5cm} 

 $=$  $(-x)$  $
	\begin{vmatrix} 
	1 & 1 &  1  \cdots  1  &  1 & 1  & 1  \cdots  &  1  \\
	-1 & (n-2)-x & -1   \cdots  -1 & -1 & 0 &-1 \cdots  & -1 \\
          -1 & -1 & (n-2)-x   \cdots  -1 & -1 & -1 &0  \cdots  & -1 \\
           \vdots  & \cdots   & \cdots         & \cdots           & \cdots & \cdots    &   \vdots \\
            0 & -1 & -1 \cdots -1& (n-2)-x &-1&-1 \cdots  & -1  \\
            -1& 0 & -1 \cdots -1& -1& (n-2)-x &-1 \cdots  & -1 \\
             \vdots  & \cdots   & \cdots              & \cdots & \cdots & \cdots   &   \vdots \\
-1 & -1 & -1 \cdots 0& -1 &-1&-1 \cdots  & (n-2)-x
	\end{vmatrix}
	$

\vspace{0.5cm} 

$=$  $(-x)$  $
	\begin{vmatrix} 
	1 & 1 &  1  \cdots  1  &  1 & 1  & 1  \cdots  &  1  \\
	0 & (n-1)-x & 0  \cdots  0 & 0 & 1 &0 \cdots  & 0 \\
          0 & 0 & (n-1)-x   \cdots  0 & 0 & 0 &1  \cdots  & 0 \\
           \vdots  & \cdots   & \cdots         & \cdots           & \cdots & \cdots    &   \vdots \\
            0 & -1 & -1 \cdots -1& (n-2)-x &-1&-1 \cdots  & -1  \\
            0& 1 & 0 \cdots 0 & 0 & (n-1)-x &0 \cdots  & 0 \\
             \vdots  & \cdots   & \cdots              & \cdots & \cdots & \cdots   &   \vdots \\
0 & 0 & 0 \cdots 1 & 0 &0 &0  \cdots  & (n-1)-x
	\end{vmatrix}
	$
  (Apply  row operations:
\newline For $ 2  \leq i  \leq \frac{n}{2} $,   $R_i \rightarrow R_1 +  R_i $
\newline For $  \frac{n}{2}+2 \leq i  \leq n $,   $R_i \rightarrow R_1 +  R_i $)
\vspace{0.5cm}

\vspace{0.5cm} 

$=$  $(-x)$ $
	\begin{vmatrix} 
	1 & 1 &  1  \cdots  1  &  1 & 1  & 1  \cdots  &  1  \\
	0 & (n-1)-x & 0  \cdots  0 & 0 & 1 &0 \cdots  & 0 \\
          0 & 0 & (n-1)-x   \cdots  0 & 0 & 0 &1  \cdots  & 0 \\
           \vdots  & \cdots   & \cdots         & \cdots           & \cdots & \cdots    &   \vdots \\
            0 & 0 & 0 \cdots 0& (n-2)-x &(n-2)-x &(n-2)-x  \cdots  & (n-2)-x   \\
            0& 1 & 0 \cdots 0 & 0 & (n-1)-x &0 \cdots  & 0 \\
             \vdots  & \cdots   & \cdots              & \cdots & \cdots & \cdots   &   \vdots \\
0 & 0 & 0 \cdots 1 & 0 &0 &0  \cdots  & (n-1)-x
	\end{vmatrix}
	$
( Apply the  row operation :  $R_{\frac{n}{2}+1} \rightarrow   R_{\frac{n}{2}+1} + ...+ R_n $)
\vspace{0.5cm} 

$=$  $(-x) ((n-2)-x)$  $
	\begin{vmatrix} 
	1 & 1 &  1  \cdots  1  &  1 & 1  & 1  \cdots  &  1  \\
	0 & (n-1)-x & 0  \cdots  0 & 0 & 1 &0 \cdots  & 0 \\
          0 & 0 & (n-1)-x   \cdots  0 & 0 & 0 &1  \cdots  & 0 \\
           \vdots  & \cdots   & \cdots         & \cdots           & \cdots & \cdots    &   \vdots \\
            0 & 0 & 0 \cdots 0& 1 &1&1\cdots  & 1   \\
            0& 1 & 0 \cdots 0 & 0 & (n-1)-x &0 \cdots  & 0 \\
             \vdots  & \cdots   & \cdots              & \cdots & \cdots & \cdots   &   \vdots \\
0 & 0 & 0 \cdots 1 & 0 &0 &0  \cdots  & (n-1)-x
	\end{vmatrix}
	$
( Apply the  row operation : $R_{\frac{n}{2}+1} \rightarrow   R_{\frac{n}{2}+1} / ((n-2)-x) $)

\vspace{0.5cm} 

$=$  $(-x) ((n-2)-x)$  $
	\begin{vmatrix} 
	1 & 1 &  1  \cdots  1  &  1 & 1  & 1  \cdots  &  1  \\
	0 & n-x & 0  \cdots  0 & 0 & n-x &0 \cdots  & 0 \\
          0 & 0 & n-x    \cdots  0 & 0 & 0 & n-x  \cdots  & 0 \\
           \vdots  & \cdots   & \cdots         & \cdots           & \cdots & \cdots    &   \vdots \\
            0 & 0 & 0 \cdots 0& 1 &1&1\cdots  & 1   \\
            0& 1 & 0 \cdots 0 & 0 & (n-1)-x &0 \cdots  & 0 \\
             \vdots  & \cdots   & \cdots              & \cdots & \cdots & \cdots   &   \vdots \\
0 & 0 & 0 \cdots 1 & 0 &0 &0  \cdots  & (n-1)-x
	\end{vmatrix}
	$

(Apply Row operation: For $ 2  \leq i  \leq \frac{n}{2} $,   $R_i \rightarrow R_i +   R_{\frac{n}{2}+i}  $)

\vspace{0.5cm} 

$=$  $(-x) ((n-2)-x) (n-x)^{{\frac{n}{2}} -1}$  $
	\begin{vmatrix} 
	1 & 1 &  1  \cdots  1  &  1 & 1  & 1  \cdots  &  1  \\
	0 & 1 & 0  \cdots  0 & 0 & 1 &0 \cdots  & 0 \\
          0 & 0 & 1   \cdots  0 & 0 & 0 & 1 \cdots  & 0 \\
           \vdots  & \cdots   & \cdots         & \cdots           & \cdots & \cdots    &   \vdots \\
            0 & 0 & 0 \cdots 0& 1 &1&1\cdots  & 1   \\
            0& 1 & 0 \cdots 0 & 0 & (n-1)-x &0 \cdots  & 0 \\
             \vdots  & \cdots   & \cdots              & \cdots & \cdots & \cdots   &   \vdots \\
0 & 0 & 0 \cdots 1 & 0 &0 &0  \cdots  & (n-1)-x
	\end{vmatrix}
	$
(Apply Row operation: For $ 2  \leq i  \leq \frac{n}{2} $,   $R_i \rightarrow R_i /(n-x) $)
\vspace{0.5cm} 

$=$  $(-x) ((n-2)-x) (n-x)^{{\frac{n}{2}} -1}$  $
	\begin{vmatrix} 
	1 & 1 &  1  \cdots  1  &  1 & 1  & 1  \cdots  &  1  \\
	0 & 1 & 0  \cdots  0 & 0 & 1 &0 \cdots  & 0 \\
          0 & 0 & 1   \cdots  0 & 0 & 0 & 1 \cdots  & 0 \\
           \vdots  & \cdots   & \cdots         & \cdots           & \cdots & \cdots    &   \vdots \\
            0 & 0 & 0 \cdots 0& 1 &1&1\cdots  & 1   \\
            0& 0 & 0 \cdots 0 & 0 & (n-2)-x &0 \cdots  & 0 \\
             \vdots  & \cdots   & \cdots              & \cdots & \cdots & \cdots   &   \vdots \\
0 & 0 & 0 \cdots 0 & 0 &0 &0  \cdots  & (n-2)-x
	\end{vmatrix}
	$
(Apply the row operation:  For $ 2  \leq i  \leq \frac{n}{2} $,  $R_{\frac{n}{2}+i} \rightarrow  R_{\frac{n}{2}+i} - R_i  $)

According to Schur's decomposition theorem  \cite {CS}, we have 
\newline $\mid L(\Delta(G)-xI) \mid$ $=$  $  (x)     ((n-2)-x)^{\frac{n}{2}}   (n-x)^{{\frac{n}{2}} -1}$ .
\newline So the  characteristic polynomial   is  $  (x)     (x-(n-2))^{\frac{n}{2}}   (x-(n))^{{\frac{n}{2}} -1}$. This completes the proof.
 \end{proof}
\begin{cor}\label{1} Let $n\geq 4$ be an even integer,  the Laplacian spectrum of $\Delta(G)$  is given by 
$ \begin{pmatrix}
0 & n-2 & n \\
1 & \frac{n}{2}  &  \frac{n}{2}-1
\end{pmatrix} $.

\end{cor}
 Based on   \cite[ Theorem 3.1]{AH}  we have the following corollary 

\begin{cor}\label{2} Let $n\geq 4$ be an even integer, the distance Laplacian spectrum of  $\Delta(G)$  is given by 
$ \begin{pmatrix}
0 & n+2 & n \\
1 & \frac{n}{2}  &  \frac{n}{2}-1
\end{pmatrix} $.

\end{cor}

According to   \cite [Theorem 2.8]{Bra}, the following corollary follows:

\begin{cor}\label{3} Let $n\geq 4$ be an even integer, the number of spanning trees of $\Delta(G)$ is  $ (n-2) ^ \frac{n}{2} (n) ^ {\frac{n}{2}-1}. $

\end{cor}

\begin{thm}\label{3.3}  Let $n\geq 4$ be an even integer and $G$ be a solvable group such that $\Delta(G)$  is  $n-2$ regular graph.  Let $\Delta(G_1)$ be the super character degree graph of  $\Delta(G)$  obtained by the addition of $n_1$ edges $( 1  \leq n_1  \leq \frac{n}{2} ) $. The characteristic polynomial of the  Laplacian matrix of $\Delta(G_1)$  is given by
\newline  $  x     (x-(n-2))^{\frac{n}{2}-n_1}   (x - n)^{({\frac{n}{2}+n_1}) -1}$ .
\end{thm}
\begin{proof} The Laplacian matrix of $\Delta(G)$ is an  $n\times n$ matrix. The rows and columns are indexed according to the vertices  $ v_1,  v_2,  ... , v_{n_1+1},   v_{n_1+2} ,...,  v_{\frac{n}{2}+1},...,
 \newline v_{\frac{n}{2}+n_1+1},  ..., v_n$. In $n-2$ regular graph,  each vertex is not adjacent with one vertex. Let us say that the vertex $v_i$ $(1\leq i \leq {\frac{n}{2}} )$ is not adjacent with the vertex  $v_{i+{\frac{n}{2}}}$. Let $n_1$  $( 1  \leq n_1  \leq \frac{n}{2} ) $ be the number of edges added to $\Delta(G)$. That is, For $1  \leq i \leq n_1 $, each $v_i$ is adjacent with $v_{\frac{n}{2}+i}$. Then

$L(\Delta(G_1))$ $= $ $ \begin{bmatrix} 
	n-1 & -1 \cdots -1& -1 \cdots  &-1 \cdots  &-1   \cdots  & -1  \\
	-1 & n-1 \cdots -1& -1 \cdots  &-1 \cdots  &-1 \cdots  & -1  \\
          \vdots  & \cdots   & \cdots   & \cdots  & \cdots       &   \vdots \\
        -1 & -1 \cdots n-2& -1 \cdots  &-1 \cdots  &0   \cdots  & -1  \\
         -1 & -1 \cdots -1& n-2 \cdots  &-1 \cdots  &-1  \cdots  & -1  \\
           \vdots  & \cdots   & \cdots   & \cdots  & \cdots       &   \vdots \\
          -1 & -1 \cdots -1& -1 \cdots  &n-1\cdots  &-1  \cdots  & -1  \\
          
              \vdots  & \cdots   & \cdots   & \cdots  & \cdots       &   \vdots \\
           -1 & -1 \cdots 0& -1 \cdots  &-1 \cdots  &n-2   \cdots  & -1  \\
           
              \vdots  & \cdots   & \cdots   & \cdots  & \cdots       &   \vdots \\
             -1 & -1 \cdots -1& -1 \cdots & -1  \cdots  &-1   \cdots  & n-2  \\
	\end{bmatrix}$
\vspace{0.5cm} 	

 So that the characteristic polynomial of  $L(\Delta(G_1))$ is  $\mid L(\Delta(G_1)-xI) \mid$

$= $ $ \begin{vmatrix} 
	(n-1)-x & -1 \cdots -1& -1 \cdots  &-1 \cdots  &-1   \cdots  & -1  \\
	-1 &(n-1)-x \cdots -1& -1 \cdots  &-1 \cdots  &-1 \cdots  & -1  \\
           \vdots  & \cdots   & \cdots   & \cdots  & \cdots       &   \vdots \\
        -1 & -1 \cdots (n-2)-x & -1 \cdots  &-1 \cdots  &0   \cdots  & -1  \\
         -1 & -1 \cdots -1& (n-2)-x \cdots  &-1 \cdots  &-1 \cdots  & -1  \\
            \vdots  & \cdots   & \cdots   & \cdots  & \cdots       &   \vdots \\
          -1 & -1 \cdots -1& -1 \cdots  &(n-1)-x  \cdots  &-1 \cdots  & -1  \\
           
            \vdots  & \cdots   & \cdots   & \cdots  & \cdots       &   \vdots \\
            -1 & -1 \cdots 0& -1 \cdots  &-1\cdots  &(n-2)-x    \cdots  & -1  \\
             
              \vdots  & \cdots   & \cdots   & \cdots  & \cdots        &   \vdots \\
             -1 & -1 \cdots -1& -1 \cdots & -1 \cdots  &-1   \cdots  & (n-2)-x  \\
	\end{vmatrix}$

\vspace{0.5cm}

$= $ $ \begin{vmatrix} 
	-x & -x \cdots -x & -x  \cdots  &-x  \cdots  &-x    \cdots  & -x   \\
	-1 &(n-1)-x \cdots -1& -1 \cdots  &-1 \cdots  &-1 \cdots  & -1  \\
           \vdots  & \cdots   & \cdots   & \cdots  & \cdots       &   \vdots \\
        -1 & -1 \cdots (n-2)-x & -1 \cdots  &-1 \cdots  &0   \cdots  & -1  \\
         -1 & -1 \cdots -1& (n-2)-x \cdots  &-1 \cdots  &-1 \cdots  & -1  \\
            \vdots  & \cdots   & \cdots   & \cdots  & \cdots       &   \vdots \\
          -1 & -1 \cdots -1& -1 \cdots  &(n-1)-x  \cdots  &-1 \cdots  & -1  \\
           
            \vdots  & \cdots   & \cdots   & \cdots  & \cdots       &   \vdots \\
            -1 & -1 \cdots 0& -1 \cdots  &-1\cdots  &(n-2)-x    \cdots  & -1  \\
             
              \vdots  & \cdots   & \cdots   & \cdots  & \cdots        &   \vdots \\
             -1 & -1 \cdots -1& -1 \cdots & -1 \cdots  &-1   \cdots  & (n-2)-x  \\
	\end{vmatrix}$
( Apply the row operation: $R_1 \rightarrow R_1 + R_2 + R_3 + ... +R_n$ )
\vspace{0.5cm}

$= $  $(-x)$  $\begin{vmatrix} 
	1& 1\cdots 1 & 1  \cdots  &1  \cdots  &1    \cdots  & 1  \\
	-1 &(n-1)-x \cdots -1& -1 \cdots  &-1 \cdots  &-1 \cdots  & -1  \\
           \vdots  & \cdots   & \cdots   & \cdots  & \cdots       &   \vdots \\
        -1 & -1 \cdots (n-2)-x & -1 \cdots  &-1 \cdots  &0   \cdots  & -1  \\
         -1 & -1 \cdots -1& (n-2)-x \cdots  &-1 \cdots  &-1 \cdots  & -1  \\
            \vdots  & \cdots   & \cdots   & \cdots  & \cdots       &   \vdots \\
          -1 & -1 \cdots -1& -1 \cdots  &(n-1)-x  \cdots  &-1 \cdots  & -1  \\
           
            \vdots  & \cdots   & \cdots   & \cdots  & \cdots       &   \vdots \\
            -1 & -1 \cdots 0& -1 \cdots  &-1\cdots  &(n-2)-x    \cdots  & -1  \\
             
              \vdots  & \cdots   & \cdots   & \cdots  & \cdots        &   \vdots \\
             -1 & -1 \cdots -1& -1 \cdots & -1 \cdots  &-1   \cdots  & (n-2)-x  \\
	\end{vmatrix}$

\vspace{0.5cm}

$= $  $(-x)$ $ \begin{vmatrix} 
	1 & 1 \cdots 1 & 1  \cdots  &1  \cdots  &1    \cdots  & 1   \\
	0 & n-x \cdots 0& 0 \cdots  &0 \cdots  &0 \cdots  & 0  \\
           \vdots  & \cdots   & \cdots   & \cdots  & \cdots       &   \vdots \\
        0 & 0 \cdots (n-1)-x & 0 \cdots  &0 \cdots  &1   \cdots  & 0  \\
         0 & 0 \cdots 0& (n-1)-x \cdots  &0 \cdots  &0 \cdots  & 0  \\
            \vdots  & \cdots   & \cdots   & \cdots  & \cdots       &   \vdots \\
          0 & 0 \cdots 0& 0 \cdots  & n-x  \cdots  &0 \cdots  & 0  \\
           
            \vdots  & \cdots   & \cdots   & \cdots  & \cdots       &   \vdots \\
            0 & 0 \cdots 1 & 0 \cdots  &0\cdots  &(n-1)-x    \cdots  & 0  \\
             
              \vdots  & \cdots   & \cdots   & \cdots  & \cdots        &   \vdots \\
             0 & 0 \cdots 0& 0 \cdots & 0 \cdots  &0   \cdots  & (n-1)-x  \\
	\end{vmatrix}$
(Apply the  row operations: For $ 2  \leq i  \leq  n $,   $R_i \rightarrow R_1 +  R_i $)

\vspace{0.5cm}

$= $  $(-x)$  $(n-x)^{(2n_1-1)}$  $ \begin{vmatrix} 
	1 & 1 \cdots 1 & 1  \cdots  &1  \cdots  &1    \cdots  & 1   \\
	0 & 1\cdots 0& 0 \cdots  &0 \cdots  &0 \cdots  & 0  \\
           \vdots  & \cdots   & \cdots   & \cdots  & \cdots       &   \vdots \\
        0 & 0      \cdots (n-1)-x & 0 \cdots  &0 \cdots  &1   \cdots  & 0  \\
         0 & 0 \cdots 0& (n-1)-x \cdots  &0 \cdots  &0 \cdots  & 0  \\
            \vdots  & \cdots   & \cdots   & \cdots  & \cdots       &   \vdots \\
           0 & 0 \cdots 0& 0 \cdots  & n-x  \cdots  &0 \cdots  & 0  \\
           
            \vdots  & \cdots   & \cdots   & \cdots  & \cdots       &   \vdots \\
            0 & 0 \cdots 1 & 0 \cdots  &0\cdots  &(n-1)-x    \cdots  & 0  \\
             
              \vdots  & \cdots   & \cdots   & \cdots  & \cdots        &   \vdots \\
             0 & 0 \cdots 0& 0 \cdots & 0 \cdots  &0   \cdots  & (n-1)-x  \\
	\end{vmatrix}$

(Apply the  row operation: 
\newline For $ 2  \leq i  \leq  n_1  $,   $R_i \rightarrow   R_i / (n-x) $)

\vspace{0.5cm}

$= $  $(-x)$  $(n-x)^{(2n_1-1)}$  $ \begin{vmatrix} 
	1 & 1 \cdots 1 & 1  \cdots  &1  \cdots  &1    \cdots  & 1   \\
	0 & 1\cdots 0& 0 \cdots  &0 \cdots  &0 \cdots  & 0  \\
           \vdots  & \cdots   & \cdots   & \cdots  & \cdots       &   \vdots \\
        0 & 0      \cdots (n-1)-x & 0 \cdots  &0 \cdots  &1   \cdots  & 0  \\
         0 & 0 \cdots 0& (n-1)-x \cdots  &0 \cdots  &0 \cdots  & 0  \\
            \vdots  & \cdots   & \cdots   & \cdots  & \cdots       &   \vdots \\
           0 & 0 \cdots 0& 0 \cdots  & 1  \cdots  &0 \cdots  & 0  \\
           
            \vdots  & \cdots   & \cdots   & \cdots  & \cdots       &   \vdots \\
            0 & 0 \cdots 1 & 0 \cdots  &0\cdots  &(n-1)-x    \cdots  & 0  \\
             
              \vdots  & \cdots   & \cdots   & \cdots  & \cdots        &   \vdots \\
             0 & 0 \cdots 0& 0 \cdots & 0 \cdots  &0   \cdots  & (n-1)-x  \\
	\end{vmatrix}$

(Apply the  row operation:
\newline For $  \frac{n}{2}+1   \leq i  \leq  \frac{n}{2}+n_1    $,   $R_i \rightarrow   R_i / (n-x) $)

\vspace{0.5cm}

$= $  $(-x)$  $(n-x)^{(2n_1-1)}$   $(n-x)^{( \frac{n}{2}-n_1)}$ $ \begin{vmatrix} 
	1 & 1 \cdots 1 & 1  \cdots  &1  \cdots  &1    \cdots  & 1   \\
	0 & 1\cdots 0& 0 \cdots  &0 \cdots  &0 \cdots  & 0  \\
           \vdots  & \cdots   & \cdots   & \cdots  & \cdots       &   \vdots \\
        0 & 0      \cdots 1 & 0 \cdots  &0 \cdots  & 1  \cdots  & 0  \\
         0 & 0 \cdots 0&  1 \cdots  &0 \cdots  &0 \cdots  & 0  \\
            \vdots  & \cdots   & \cdots   & \cdots  & \cdots       &   \vdots \\
           0 & 0 \cdots 0& 0 \cdots  & 1 \cdots  &0 \cdots  & 0  \\
           
            \vdots  & \cdots   & \cdots   & \cdots  & \cdots       &   \vdots \\
            0 & 0 \cdots 1 & 0 \cdots  &0\cdots  &(n-1)-x    \cdots  & 0  \\
             
              \vdots  & \cdots   & \cdots   & \cdots  & \cdots        &   \vdots \\
             0 & 0 \cdots 0& 0 \cdots & 0 \cdots  &0   \cdots  & (n-1)-x  \\
	\end{vmatrix}$

(Apply the  row operations consecutively, 
\newline For $ n_1+1 \leq i  \leq   \frac{n}{2}  $,   $R_i \rightarrow R_i +  R_{ \frac{n}{2}+i} $
\newline For $ n_1+1 \leq i  \leq   \frac{n}{2}  $,   $R_i \rightarrow   R_i / (n-x) $)

\vspace{0.5cm}

$= $  $(-x)$ $ (n-x)^{({\frac{n}{2}+n_1}) -1}$  $ \begin{vmatrix} 
	1 & 1 \cdots 1 & 1  \cdots  &1  \cdots  &1    \cdots  & 1   \\
	0 & 1\cdots 0& 0 \cdots  &0 \cdots  &0 \cdots  & 0  \\
           \vdots  & \cdots   & \cdots   & \cdots  & \cdots       &   \vdots \\
        0 & 0      \cdots 1 & 0 \cdots  &0 \cdots  & 1  \cdots  & 0  \\
         0 & 0 \cdots 0&  1 \cdots  &0 \cdots  &0 \cdots  & 0  \\
            \vdots  & \cdots   & \cdots   & \cdots  & \cdots       &   \vdots \\
           0 & 0 \cdots 0& 0 \cdots  & 1 \cdots  &0 \cdots  & 0  \\
           
            \vdots  & \cdots   & \cdots   & \cdots  & \cdots       &   \vdots \\
            0 & 0 \cdots 1 & 0 \cdots  &0\cdots  &(n-1)-x    \cdots  & 0  \\
             
              \vdots  & \cdots   & \cdots   & \cdots  & \cdots        &   \vdots \\
             0 & 0 \cdots 0& 0 \cdots & 0 \cdots  &0   \cdots  & (n-1)-x  \\
	\end{vmatrix}$

\vspace{0.5cm}

$= $  $(-x)$ $ (n-x)^{({\frac{n}{2}+n_1}) -1}$  $ \begin{vmatrix} 
	1 & 1 \cdots 1 & 1  \cdots  &1  \cdots  &1    \cdots  & 1   \\
	0 & 1\cdots 0& 0 \cdots  &0 \cdots  &0 \cdots  & 0  \\
           \vdots  & \cdots   & \cdots   & \cdots  & \cdots       &   \vdots \\
        0 & 0      \cdots 1 & 0 \cdots  &0 \cdots  & 1  \cdots  & 0  \\
         0 & 0 \cdots 0&  1 \cdots  &0 \cdots  &0 \cdots  & 0  \\
            \vdots  & \cdots   & \cdots   & \cdots  & \cdots       &   \vdots \\
           0 & 0 \cdots 0& 0 \cdots  & 1 \cdots  &0 \cdots  & 0  \\
           
            \vdots  & \cdots   & \cdots   & \cdots  & \cdots       &   \vdots \\
            0 & 0 \cdots 0 & 0 \cdots  &0\cdots  &(n-2)-x    \cdots  & 0  \\
             
              \vdots  & \cdots   & \cdots   & \cdots  & \cdots        &   \vdots \\
             0 & 0 \cdots 0& 0 \cdots & 0 \cdots  &0   \cdots  & (n-2)-x  \\
	\end{vmatrix}$
(Apply the row operation:  For $ n_1+1 \leq i  \leq   \frac{n}{2}  $,   $ R_{ \frac{n}{2}+i} \rightarrow   R_{ \frac{n}{2}+i} - R_i  $)

By using Schur’s decomposition theorem \cite {CS}, we have 
\newline $\mid L(\Delta(G_1)-xI)\mid$ $=$  $(-x) ((n-2)-x)^{\frac{n}{2}-n_1}  (n-x)^{({\frac{n}{2}+n_1}) -1}$.
\newline So the  characteristic polynomial   is  $  x     (x-(n-2))^{\frac{n}{2}-n_1}   (x-(n))^{({\frac{n}{2}+n_1}) -1}$   and hence the result holds.
\end {proof}
\begin{cor}\label{3}Let  $n\geq 4$ be an even integer,  the Laplacian spectrum of $\Delta(G_1)$  is given by 
$ \begin{pmatrix}
0 & n-2 & n \\
1 & \frac{n}{2} - n_1  &  (\frac{n}{2}+n_1)-1
\end{pmatrix} $.

\end{cor}

\cite[ Theorem 3.1]{AH} demonstrates the following corollary:

\begin{cor}\label{4}Let  $n\geq 4$ be an even integer,  the distance Laplacian spectrum of $\Delta(G_1)$   is given by 
$ \begin{pmatrix}
0 & n + 2 & n \\
1 & \frac{n}{2} - n_1  &  (\frac{n}{2}+n_1)-1
\end{pmatrix} $.

\end{cor}

  \cite [Theorem 2.8]{Bra} implies following corollary:
\begin{cor}\label{5} Let $n\geq 4$ be an even integer, the number of spanning trees of $\Delta(G_1)$ is  $ (n-2) ^ {(\frac{n}{2} - n_1)}  (n) ^ {(\frac{n}{2} + n_1) - 1}. $

\end{cor}

\begin{thm}\label{3.4}  Let $n\geq 3$ be an  integer and $G$ be a solvable group. Let $\Delta(G)$ be a character degree graph with diameter $2$ and it is not a block. $\Delta(G)$  has the structure:  $\Delta(G) = K_{n-(n_1+1)} - v - K_{n_1}$, for some integer $1\leq n_1\leq n-(n_1+1)$, where $V( K_{n-(n_1+1)}) \cup \{v\}$ and $V( K_{n_1}) \cup \{v\}$ generates a clique. The characteristic polynomial of $L (\Delta(G)) $ is given by 
\newline $ (x) (x-1) (x-n) (x-(n_1+1))^{n_1-1} (x-(n-n_1))^{n-(n_1+2)}$.
\end{thm}
\begin{proof} The Laplacian matrix of $\Delta(G)$ is  $n\times n$ matrix. The rows and columns are indexed in order by the vertices $ v_1,  v_2, ... , v_{n_1},  v_{n_1+1},  v_{n_1+2},  . . . , v_n$. Complete graph   $K_{n_1}$ is not adjacent with any vertex of $ K_{n-(n_1+1)}$ and vice versa. 

\vspace{0.5cm}

$L(\Delta(G))$ $= $ $ \begin{bmatrix} 
	n_1 & -1   \cdots -1& -1 & 0\cdots  & 0   \\
       -1 &n_1   \cdots -1& -1 & 0 \cdots  & 0   \\

          \vdots  & \cdots   & \cdots           &   \vdots \\

         -1 & -1   \cdots n_1 & -1 & 0 \cdots  & 0   \\
           -1 & -1   \cdots -1 & n-1 &  -1  \cdots  & -1   \\
            0 & 0  \cdots     0 & -1 & n-(n_1+1) \cdots  & -1   \\
            
              \vdots  & \cdots      & \cdots  & \cdots       &   \vdots \\

              0 & 0   \cdots     0 & -1 & -1  \cdots  & n-(n_1+1)   \\

	\end{bmatrix}$
\vspace{0.5cm} 	

 The characteristic polynomial of $L(\Delta(G))$ is $\mid L(\Delta(G)-xI) \mid$

$= $ $ \begin{vmatrix} 

	n_1-x & -1   \cdots -1& -1 & 0 \cdots  & 0   \\
       -1 &n_1-x  \cdots -1& -1 & 0 \cdots  & 0   \\

          \vdots  & \cdots   & \cdots            &   \vdots \\

         -1 & -1  \cdots n_1-x & -1 & 0 \cdots  & 0   \\
           -1 & -1   \cdots -1 & (n-1)-x &  -1  \cdots  & -1   \\
            0 & 0   \cdots     0 & -1 & (n-(n_1+1))-x  \cdots  & -1   \\
            
              \vdots  & \cdots      & \cdots  & \cdots       &   \vdots \\

              0 & 0  \cdots     0 & -1 & -1  \cdots  & (n-(n_1+1))-x   \\

	\end{vmatrix}$

\vspace{0.5cm} 	

$= $ $ \begin{vmatrix} 

	-x & -x    \cdots -x  & -x  &-x  \cdots  & -x   \\
       -1 &n_1-x  \cdots -1& -1 & 0 \cdots  & 0   \\

          \vdots  & \cdots   & \cdots            &   \vdots \\

         -1 & -1  \cdots n_1-x & -1 & 0 \cdots  & 0   \\
           -1 & -1   \cdots -1 & (n-1)-x &  -1  \cdots  & -1   \\
            0 & 0   \cdots     0 & -1 & (n-(n_1+1))-x  \cdots  & -1   \\
            
              \vdots  & \cdots      & \cdots  & \cdots       &   \vdots \\

              0 & 0  \cdots     0 & -1 & -1  \cdots  & (n-(n_1+1))-x   \\

	\end{vmatrix}$
(Apply the row operation: $R_1 \rightarrow R_1 + R_2 + R_3 + ... +R_n$ )
\vspace{0.5cm}

$= $  $(-x)$  $\begin{vmatrix} 

	1 & 1    \cdots 1  & 1  & 1 \cdots  & 1  \\
       -1 &n_1-x  \cdots -1& -1 & 0 \cdots  & 0   \\

          \vdots  & \cdots   & \cdots            &   \vdots \\

         -1 & -1  \cdots n_1-x & -1 & 0 \cdots  & 0   \\
           -1 & -1   \cdots -1 & (n-1)-x &  -1  \cdots  & -1   \\
            0 & 0   \cdots     0 & -1 & (n-(n_1+1))-x  \cdots  & -1   \\
            
              \vdots  & \cdots      & \cdots  & \cdots       &   \vdots \\

              0 & 0  \cdots     0 & -1 & -1  \cdots  & (n-(n_1+1))-x   \\
	\end{vmatrix}$

\vspace{0.5cm}

$= $  $(-x) (n-x) $  $\begin{vmatrix} 

	 1 & 1    \cdots 1  & 1  & 1 \cdots  & 1  \\
         0 &(n_1+1)-x  \cdots 0& 0 & 1 \cdots  & 1   \\

          \vdots  & \cdots   & \cdots            &   \vdots \\

         0 & 0  \cdots (n_1+1)-x & 0 & 1 \cdots  & 1  \\
           0 & 0   \cdots 0 & 1 &  0 \cdots  & 0  \\
            0 & 0   \cdots     0 & -1 & (n-(n_1+1))-x  \cdots  & -1   \\
            
              \vdots  & \cdots       & \cdots       &   \vdots \\

              0 & 0  \cdots     0 & -1 & -1  \cdots  & (n-(n_1+1))-x   \\
            \end{vmatrix}$
\vspace{0.2cm}
 ( Apply the  row operations consecutively, 
\newline For $ 2  \leq i  \leq  n_1+1 $,   $R_i \rightarrow R_1 +  R_i $
\newline  $R_{ n_1+1}\rightarrow  R_{ n_1+1} / (n-x) $)

\vspace{0.5cm}

$= $  $(-x) (n-x) $  $\begin{vmatrix} 

	 1 & 1    \cdots 1  & 1  & 1 \cdots  & 1  \\
         0 &(n_1+1)-x  \cdots 0& 0 & 1 \cdots  & 1   \\

          \vdots  & \cdots   & \cdots            &   \vdots \\

         0 & 0  \cdots (n_1+1)-x & 0 & 1 \cdots  & 1  \\
           0 & 0   \cdots 0 & 1 &  0 \cdots  & 0  \\
            0 & 0   \cdots     0 & 0 & (n-(n_1+1))-x  \cdots  & -1   \\
            
              \vdots  & \cdots       & \cdots       &   \vdots \\

              0 & 0  \cdots     0 & 0 & -1  \cdots  & (n-(n_1+1))-x   \\
            \end{vmatrix}$
\vspace{1cm}

( Apply the  row operation:  For $ n_1+2  \leq i  \leq  n $,   $R_i \rightarrow R_{n_1+1} +  R_i $)

\vspace{0.5cm}

$= $  $(-x) (n-x) (1-x) $  $\begin{vmatrix} 

	 1 & 1    \cdots 1  & 1  & 1 \cdots  & 1  \\
         0 &(n_1+1)-x  \cdots 0& 0 & 1 \cdots  & 1   \\

          \vdots  & \cdots   & \cdots            &   \vdots \\

         0 & 0  \cdots (n_1+1)-x & 0 & 1 \cdots  & 1  \\
           0 & 0   \cdots 0 & 1 &  0 \cdots  & 0  \\
            0 & 0   \cdots     0 & 0 & 1  \cdots  & 1  \\
            
              \vdots  & \cdots       & \cdots       &   \vdots \\

              0 & 0  \cdots     0 & 0 & -1  \cdots  & (n-(n_1+1))-x   \\
            \end{vmatrix}$

( Apply the  row operations consecutively, 
\newline   $R_{n_1+2} \rightarrow R_{n_1+2} +  R_{n_1+3}+...+R_n $
\newline   $R_{n_1+2} \rightarrow R_{n_1+2} / (1-x)$)

$= $  $(-x) (n-x) (1-x) $  $\begin{vmatrix} 

	 1 & 1    \cdots 1  & 1  & 1 \cdots  & 1  \\
         0 &(n_1+1)-x  \cdots 0& 0 & 1 \cdots  & 1   \\

          \vdots  & \cdots   & \cdots            &   \vdots \\

         0 & 0  \cdots (n_1+1)-x & 0 & 1 \cdots  & 1  \\
           0 & 0   \cdots 0 & 1 &  0 \cdots  & 0  \\
            0 & 0   \cdots     0 & 0 & 1  \cdots  & 1  \\
            
              \vdots  & \cdots       & \cdots       &   \vdots \\

              0 & 0  \cdots     0 & 0 & 0  \cdots  & (n-n_1)-x   \\
            \end{vmatrix}$

\vspace{0.75cm}( Apply the  row operation:  For $ n_1+2  \leq i  \leq  n $,   $R_i \rightarrow R_i+ R_{n_1+2}$ )

By using Schur’s decomposition theorem \cite {CS}, we have 
\newline $\mid L(\Delta(G)-xI)\mid $ $=$ $(-x) (n-x) (1-x) ((n_1+1)-x)^{n_1-1} ((n-n_1)-x)^{n-(n_1+2)}. $ 
\newline So the  characteristic polynomial   is
 \newline $ (x) (x-1) (x-n) (x-(n_1+1))^{n_1-1} (x-(n-n_1))^{n-(n_1+2)}$  and hence the result holds.

\begin{cor}\label{5}Let  $n\geq 3$ be an  integer,  the Laplacian spectrum of $\Delta(G)$  is given by 
$ \begin{pmatrix}
0 & 1& n & n_1+1 & n - n_1\\
1 & 1 & 1 & n_1- 1 & n - (n_1+2)
\end{pmatrix} $.

\end{cor}

 According to  \cite[ Theorem 3.1]{AH}  the following corollary follows:

\begin{cor}\label{6}Let  $n\geq 3$ be an  integer,  the distance Laplacian spectrum of $\Delta(G)$  is given by 
$ \begin{pmatrix}
0 & 2n-1& n & 2n - (n_1+1) & n + n_1\\
1 & 1 & 1 & n_1- 1 & n - (n_1+2)
\end{pmatrix} $.

\end{cor}
From  \cite [Theorem 2.8]{Bra} follows the following corollary:
\begin{cor}\label{7} Let  $n\geq 3$ be an  integer,  the number of spanning trees of $\Delta(G)$ is  $ n  (n_1 + 1) ^ {( n_1 - 1)}  (n - n_1) ^ {(n-(n_1 + 2))} $.

\end{cor}

\end{proof}

\end{document}